\documentclass{amsart}
\usepackage{amsfonts,amssymb,amsmath,amsthm}
\usepackage{url}
\usepackage{enumerate}

\newtheorem{theorem}{Theorem}[section]

\newcommand{\q}{\quad}

\newcommand{\wh}{\widehat}

\newcommand{\tf}{\tfrac}

\author{Geoff Diestel}
\address{
Texas A\&M University-Central Texas\\
1001 Leadership Place, Killeen TX 76549, USA}
\email{gdiestel@tamuct.edu}

\begin{document}

\title[Running Head]{An extension of Nikishin's factorization theorem}

\begin{abstract}
A Nikishin-Maurey characterization is given for bounded subsets of weak-type Lebesgue spaces. New factorizations for linear and multilinear operators are shown to follow.

\end{abstract}
\maketitle


\section{Introduction}

Let $(\Omega,\mu)$ be a non-atomic probability space. The space of scalar-valued $\mu$-measurable functions $L_0=L_0(d\mu)$ is equipped with the topology of convergence in measure and $f=g$ in $L_0$ if $f(\omega)=g(\omega)$ for $\mu$-almost every $\omega \in \Omega$.

For $0<p<\infty$, the Lebesgue space $L_p$ is defined by
$$f\in L_p \Leftrightarrow \|f\|_p=\Big(\int|f|^p\ d\mu\Big)^{1/p}<\infty.$$
Of course $\|\cdot\|_p$ is a quasi-norm making $L_p$ into a quasi-Banach space and $L_p$ is a Banach space if $1\le p\le \infty$ where $L_{\infty}$ is the space of essentially bounded functions equipped with the essential supremum norm.

For $0<p<\infty$, the weak Lebesgue space $L_{p,\infty}$ is defined by
$$f\in L_{p,\infty}\Leftrightarrow \|f\|_{p,\infty}:=\sup_{t>0}t\mu(\{|f|>t\})^{1/p}<\infty.$$
The quantity $\|\cdot\|_{p,\infty}$ is a (complete) quasi-norm (see \cite{Grafakosbook1}, Ch 1.1.1). 

The intent of this article is to extend the existing Nikishin-Maurey theory. Although this theory is often presented in the context of factoring operators, the foundational results of this theory are concerned with bounded subsets of non-negative measurable functions.
If $\mathcal{F}$ is such a set of measurable functions and $0<p<\infty$, Nikishin \cite{Nikishin} characterized the existence of a positive measurable function $g$ such that 
$$\sup_{f\in \mathcal{F}}\|f/g\|_{p,\infty}<\infty.$$
The existence of such a function is equivalent to the existence of a decreasing function $C:(0,\infty)\rightarrow (0,\infty)$ such that $\lim_{t\rightarrow \infty} C(t)=0$ and
$$\mu(\{\sup_j|c_jf_j|>t\})\le C(t)$$
for all $t>0$, finitely supported sequences $(f_j)$ from $\mathcal{F}$ and scalars $(c_j)$ satisfying $\sum_j|c_j|^p\le 1$. In other words, $g$ exists if and only if
$$\big\{\sup_j|c_jf_j|:n\in \mathbb N, f_1,\cdots,f_n\in\mathcal{F},|c_1|^p+\cdots+|c_n|^p\le 1\big\}$$
is bounded in $L_0$. The key step in proving $g$'s existence from  the boundedness of this set is purely constructive. A nice presentation of the proof of Nikishin's \cite{Nikishin} characterization may be found in Proposition III.H.2 of \cite{Wojtaszczyk}. To summarize, for $0<\epsilon<1$, the boundedness of the above set within $L_0$ implies the existence of a subset $E_{\epsilon}$ of $\Omega$ and a constant $C_{\epsilon}<\infty$ such that $\mu(E_{\epsilon})\ge 1-\epsilon$ and
$$\sup_{f\in \mathcal{F}}\|1_{E_{\epsilon}}f\|_{p,\infty}\le C_{\epsilon}.$$
Of course the constant $C_{\epsilon}$ is tied to the function $C$ from the above maximal estimate and thus it may be that $\lim_{\epsilon\rightarrow 0}C_{\epsilon}=\infty$. If this limit were finite there would be nothing to prove as the set $\mathcal{F}$ would be bounded in $L_{p,\infty}$. 
Given any decreasing null sequence $(\epsilon_n)$ such that $0<\epsilon_n<1$, $g$ is constructed by selecting an unbounded increasing sequence of positive scalars $(a_n)$ such that $g=a_11_{E_{\epsilon_1}}+\sum_{n>1}a_n1_{E_{\epsilon_n}\setminus E_{\epsilon_{n-1}}}$
where the scalars $(a_n)$ are chosen to  counteract the growth of constants $(C_{\epsilon_n})$ in order to obtain the factorization $\mathcal{F}=g(g^{-1}\mathcal{F})$ where $g^{-1}\mathcal{F}$ is bounded in $L_{p,\infty}$. Of course, any selection of scalars $(a_n)$ defines $g$ as a positive element of $L_0$. However, it is of interest to identify a more specific space for $g$ based on the particular function $C(t)$ from the above
maximal estimate. Theorem \ref{weakN} is a formulation of Nikishin's theorem for the specific case that $\mathcal{F}$ is a subset of $L_{q,\infty}$ for some $0<q<p$. In this case, it is natural to consider
$C(t)=t^{-q}$. Much of the proof follows from a technical adaptations of the arguments used to prove Nikishin's \cite{Nikishin} theorem as presented in \cite{Wojtaszczyk}. However, the aforementioned construction of the function $g$ is insufficient to obtain $g$ from the natural weak-Lebesgue space and an extra compactness argument will be used to prove $N3\Rightarrow N1$.

\begin{theorem}\label{weakN}
Let $0<q<p<\infty$, $1/r=1/q-1/p$ and $\mathcal{F}$ be a subset of non-negative elements of $L_{q,\infty}$. Then the following conditions are equivalent.
\begin{enumerate}
\item[$N1.$] There exists a constant $C<\infty$ and a positive $g\in L_{r,\infty}$ so that $\|g\|_{r,\infty}= 1$ and
$$\sup_{f\in \mathcal{F}}\|f/g\|_{p,\infty}\le C.$$
\item[$N2.$] There exists a constant $C<\infty$ so that
$$\|\sup_j|c_jf_j|\|_{q,\infty}\le C\Big(\sum_j|c_j|^p\Big)^{1/p}$$
for all finitely supported sequences $(f_j)$ from $\mathcal{F}$ and scalars $(c_j)$.
\item[$N3.$] There exists a constant $C<\infty$ such that for any $0<\epsilon<1$ there exists a measurable set $E_{\epsilon}$ such that $\mu(E_{\epsilon})\ge 1-\epsilon$ and
$$\sup_{f\in \mathcal{F}}\|1_{E_{\epsilon}}f\|_{p,\infty}\le C\epsilon^{-1/r}.$$
\end{enumerate}
\end{theorem}

Theorem \ref{weakN} will be proven in Section 3. Section 2 contains applications related to the factorization of linear and multilinear operators. For more on Nikishin-Maurey theory and its applications to Banach space theory, the reader is referred to \cite{AlbiacKalton}, \cite{Wojtaszczyk}, \cite{handbook1}, and \cite{handbook2}.

\section{Factoring Operators}

The topology of any locally bounded topological vector space $X$ is induced by a function $\|\cdot\|:X\rightarrow [0,\infty)$ satisfying the following conditions.
\begin{enumerate}
\item $\|x\|>0$ for all $x\neq 0$.
\item $\|cx\|=|c|\|x\|$ for all scalars $c$ and all $x\in X$.
\item There exists $1\le C<\infty$ so that $\|x+y\|\le C(\|x\|+\|y\|)$ for all $x,y\in X$.
\end{enumerate}
The function $\|\cdot\|$ is called a quasi-norm and $X$ is a quasi-Banach space if $X$ is complete with respect to $\|\cdot\|$. Furthermore, if $C=1$ then $X$ is a Banach space.

For $0<p\le 2$, a quasi-Banach space $X$ has Rademacher type $p$ if there exists a constant $T_p(X)<\infty$ such that
$$\mathbb E\Big\|\sum_j\epsilon_jx_j\Big\|^p\le T_p(X)^p\sum_j\|x_j\|^p$$
for all finitely supported sequences $(x_j)$ from $X$. Here $(\epsilon_j)$ is a sequence of independent Bernoulli random variables satisfying
$$P(\epsilon_j=1)=1/2=P(\epsilon_j=-1).$$

For $k\ge 1$ and quasi-Banach space(s) $X_1,\cdots,X_k$, $T:X_1\times \cdots\times X_k\rightarrow L_0$ is $k$-sublinear if each coordinate map is sublinear. If each coordinate map is linear then $T$ is a $k$-linear operator. 
In \cite{Diestel} it was shown that for a bounded $k$-sublinear operator $T:X_1\times\cdots\times X_k\rightarrow L_0$ there exists a decreasing function $C:(0,\infty)\rightarrow (0,\infty)$ such that $\lim_{t\rightarrow \infty}C(t)=0$ and
$$\mu(\{\sup_j|T(x_{1,j},\cdots,x_{k,j})|>t\})\le C(t)$$
for all finitely supported sequences $(x_{i,j})_j$ from $X_i$, $1\le i\le k$, such that
$$\sum_j(\|x_{1,j}\|\cdots\|x_{k,j}\|)^p\le 1$$
where $1/p=1/p_1+\cdots+1/p_k$ and $0<p_1,\cdots,p_k\le 2$ are the respective Rademacher types of the quasi-Banach spaces $X_1,\cdots,X_k$. Thus, by Nikishin's \cite{Nikishin} characterization there exists a positive measurable $g$ so that
$$\|g^{-1}T(x_1,\cdots,x_k)\|_{p,\infty}\le \|x_1\|\cdots\|x_k\|$$
for all $x_i\in X_i$, $1\le i\le k$. Of course the case $k=1$ is due to Nikishin \cite{Nikishin}. Analogous results are also shown in \cite{Diestel} for operators mapping into $L_q$ for some $0<q<p$. These results utilize Pisier's \cite{Pisier} characterization for the factorization of subsets of $L_q$ through $L_{p,\infty}$ if $0<q<p$. 

Maurey \cite{Maurey} characterized when a subset of $L_q$ may be factored through $L_p$ for $0<q<p<\infty$. Due to the stable laws and the Kahane-Khintchine inequalities, this result only applies to linear operators defined on a space of type $2$, i.e. every linear operator $T:X\rightarrow L_q$ factors through $L_2$ if $0<q<2$. There are no direct applications of Maurey's \cite{Maurey} results for factoring multilinear operators through $L_p$ using Rademacher type. However, partial analogs of Maurey's \cite{Maurey} result are shown to hold in \cite{Diestel} and \cite{Kalton2}. Suppose $X_1$ and $X_2$ have type 2 and $T:X_1\times X_2\rightarrow L_0$ is a bounded bilinear operator. Then not only does $T$ factor through $L_{1,\infty}$ but the range of $T$ is locally convex in $L_0$. This is shown in \cite{Diestel} as an application of the characterization of the Rademacher decoupling property for quasi-Banach spaces given in \cite{Kalton2}. Thus, the $T$-induced linear map is continuous from the projective tensor-product $X_1\wh{\otimes}X_2$ into $L_0$. This means that $T$ factors through a Banach space, i.e. there exists a Banach Space $Z$, a bilinear operator $B:X_1\times X_2\rightarrow Z$ and a linear operator $L:Z\rightarrow L_0$ such that $T=LB$.

As with previously existing Nikishin-Maurey theorems, Theorem \ref{weakN} applies to linear and multilinear operators due to their homogeneity (see \cite{Diestel}). Suppose $k\ge 1$ and $0<p_1,\cdots,p_k\le 2$ are the respective Rademacher types of the quasi-Banach spaces $X_1,\cdots, X_k$. If $1/p=1/p_1+\cdots+1/p_k$ and $T:X_1\times \cdots\times X_k\rightarrow L_{q,\infty}$ is a bounded $k$-linear operator and $0<q<p$ then the Kahane-Khintchine inequalities imply that there exists a constant $C=C(q,p_1,\cdots,p_k,T)<\infty$ such that
$$\Big\|\Big(\sum_{j_1,\cdots,j_k}|T(x_{1,j_1},\cdots,x_{k,j_k})|^2\Big)^{1/2}\Big\|_{q,\infty}\le C\prod_{i=1}^k\Big(\sum_{j_i}\|x_{i,j_i}\|^{p_i}\Big)^{1/p_i}$$
for all finitely supported sequences $(x_{i,j})_j$ from $X_i$, $1\le i\le k$. Of course this estimate implies that
$$\|\sup_j|T(x_{1,j},\cdots,x_{k,j})|\|_{q,\infty}\le C\prod_{i=1}^k\Big(\sum_j\|x_{i,j}\|^{p_i}\Big)^{1/p_i}$$
for all finitely supported sequences $(x_{i,j})_j$ from $X_i$, $1\le i\le k$. As illustrated in \cite{Kalton2} and \cite{Diestel}, the homogeneity of $T$ and of the identity $1/p=1/p_1+\cdots+1/p_k$ implies
$$\|\sup_j|T(x_{1,j},\cdots,x_{k,j})|\|_{q,\infty}\le C\Big(\sum_j(\|x_{1,j}\|\cdots\|x_{k,j}\|)^p\Big)^{1/p}$$
for all finitely supported sequences $(x_{i,j})_j$ from $X_i$, $1\le i\le k$. An immediate consequence of this estimate and Theorem \ref{weakN} is the following.

\begin{theorem}\label{weakNA}
Let $k\ge 1$, $0<p_1,\cdots,p_k\le 2$, $1/p=1/p_1+\cdots+1/p_k<1/q<\infty$, $1/r=1/q-1/p$ and suppose $T:X_1\times \cdots\times X_k\rightarrow L_{q,\infty}$ is a continuous $k$-linear operator where $p_1,\cdots,p_k$ are the respective Rademacher types of the quasi-Banach spaces $X_1,\cdots,X_k$.  Then there exists a positive $g\in L_{r,\infty}$ such that
$$\|g^{-1}T(x_1,\cdots,x_k)\|_{p,\infty}\le \|x_1\|\cdots\|x_k\|$$
for all $x_i\in X_i$, $1\le i\le k$.
\end{theorem}

Theorem \ref{weakNA} identifies new factorizations not established by the characterizations of Nikishin \cite{Nikishin}, Maurey \cite{Maurey}, or Pisier \cite{Pisier} due the the fact that the operator takes values in $L_{q,\infty}$ and that the factorization $T=g(g^{-1}T)$ is defined by an element of $g$ of $L_{r,\infty}$ where $1/r=1/q-1/p$.

\section{Proof of Theorem \ref{weakN}}

All Nikishin-Maurey characterizations like Theorem \ref{weakN} may be reduced to any particular value of $p$ by considering
$$\mathcal{F}^t=\{f^t:f\in \mathcal{F}\}$$
for any $t>0$. Since $1/r=1/q-1/p$, it follows that $1/(r/t)=1/(q/t)+1/(p/t)$ and all the conditions $N1-N3$ translate into equivalent conditions about $\mathcal{F}^t$. Picking $t$ so that $p/t$ has a particular value and proving the equivalence of these new conditions for $\mathcal{F}^t$ implies the general equivalence of $N1-N3$.

The following proof of $N1\Rightarrow N2\Rightarrow N3$ follows from a technical reworking of the arguments used to prove the analogous implications of Nikishin's \cite{Nikishin} theorem as presented in \cite{Wojtaszczyk}. However, the aforementioned constructive arguments used to prove the final implication of Nikshin's \cite{Nikishin} theorem are insufficient to construct the function $g$ from $N1$ because of the added requirement that $g\in L_{r,\infty}$ where $1/r=1/q-1/p$. However, assuming $N3$ and using the aforementioned reduction with the assumption that $r>2$, one may construct a sequence $(g_n)$ such that
$$\sup_n\sup_{f\in \mathcal{F}}\|f/g_n\|_{p,\infty}<\infty$$
and
$$\sup_n\|g_n\|_{r_n,\infty}<\infty$$
where $2<r_n<r$ and $(r_n)$ is increasing to $r$. Thus, $(g_n)$ is bounded in $L_2$ and the Banach-Saks Theorem guarantees that $(g_n)$ has a subsequence with convergent Ces\`aro means. The above two conditions will be used to show that the limit $g$ of these Ces\`aro means will satisfy $N1$ and this will complete the proof by establishing $N3\Rightarrow N1$.

\begin{proof}
Assuming $N1$, let $(f_j)$ be a finitely supported sequence from $\mathcal{F}$, $(c_j)$ be scalars, and define $F=\sup_j|c_jf_j|$. Then $1/r=1/q-1/p$ implies
$$\|F\|_{q,\infty}\le \|F/g\|_{p,\infty}\|g\|_{r,\infty}.$$
Since $(f_j)$ is a finitely supported sequence from $\mathcal{F}$, there is a partition $(E_j)$ of $\Omega$ so that
$$F=|c_kf_k|$$
almost everywhere on $E_k$. The disjointness of the partition implies
$$\mu(\{F/g>t\})\le \sum_k\mu(\{|c_kf_k|/g>t\}).$$
Hence, $N1$ implies
$$\mu(\{\sup_j|c_jf_j|/g>t\})\le \sum_k\mu(\{|c_kf_k|>t\})\le Ct^{-p}\sum_k|c_k|^p$$
and $N2$ follows with same constant as in $N1$ because $\|g\|_{r,\infty}=1$.

Now assume $N2$ holds. Without loss of generality, assume $C=1$ by normalizing $\mathcal{F}$.
Fix $0<\epsilon<1$ and consider a measurable set $B$ to be $\epsilon$-bad if there exists $f\in \mathcal{F}$ so that
$$\mu(B)|f(\omega)|^p>\epsilon^{-p/q}$$
for all $\omega\in B$. 

If there are no $\epsilon$-bad sets for some $0<\epsilon<1$ then for every $f\in \mathcal{F}$ and every $t>0$, there exists $\omega\in \{|f|>t\}$ so that
$$\mu(\{|f|>t\})t^p<\mu(\{|f|>t\})|f(\omega)|^p\le \epsilon^{-p/q}.$$
In this case $\mathcal{F}$ is bounded in $L_{p,\infty}$ and $N3$ follows in a trivial manner.

Suppose there are $\epsilon$-bad sets for every $0<\epsilon<1$. Fix $\epsilon$ and suppose $(B_j)$ is a maximal family of disjoint $\epsilon$-bad sets. Thus, for each $j$ there exists $f_j\in \mathcal{F}$ so that
$$\mu(B_j)|f_j(\omega)|^p>\epsilon^{-p/q}$$
for all $\omega\in B_j$. For each $j$ let $c_j=\mu(B_j)^{1/p}$. Notice that 
$$\sup_j|c_jf_j|^p>\epsilon^{-p/q}$$
everywhere on $B=\cup_jB_j$. Thus, for any $n\in \mathbb N$, $N2$ implies that
$$\mu(\cup_{j=1}^nB_j)\le \mu(\{\sup_{j\le n}|c_jf_j|>\epsilon^{-1/q}\})\le \epsilon \Big(\sum_{j=1}^n|c_j|^p\Big)^{q/p}\le \epsilon\mu(B)^{q/p}.$$
Letting $n$ tend to infinity implies that
$$\mu(B)^{1-q/p}\le \epsilon.$$
Therefore,
$$\mu(B)\le \epsilon^{p/(p-q)}.$$
So, if $E_{\epsilon^{p/(p-q)}}=\Omega\setminus B$ then $\mu(E_{\epsilon^{p/(p-q)}})\ge 1-\epsilon^{p/(p-q)}$ and $E_{\epsilon^{p/(p-q)}}$ is not $\epsilon$-bad. Therefore,
$$\|1_{E_{\epsilon^{p/(p-q)}}}f\|_{p,\infty}^p\le \epsilon^{-p/q}.$$
This holds for all $0<\epsilon<1$. By making the substitution $\delta=\epsilon^{p/(p-q)}$ then for every $0<\delta<1$, the condition $1/r=1/q-1/p$ implies that there exists $E_{\delta}$ such that $\mu(E_{\delta})\ge 1-\delta$ and
$$\|1_{E_{\delta}}f\|_{p,\infty}^p\le \delta^{-(p-q)/q}=\delta^{-p/r}.$$
By undoing the assumed normalization of $N2$, it follows that $N3$ holds with the same constant as in $N2$.

Assume $N3$ with $C=1$. Moreover, the fact that $\mathcal{F}$ is a set of positive measurable functions, the homogeneity of the desired estimates on the indices $q,p,$ and $r$ due to the identity $1/r=1/q-1/p$ allows for the assumption that $r>2$.

Fix $n$ and let $(\epsilon_{n,m})_m$ be the sequence defined by $\epsilon_{n,m}=\big(1/(2m^n)\big)_m$. For each $m$, let $E_{n,m}$ be the set $E_{\epsilon_{n,m}}$ from condition $N3$ for $\epsilon=\epsilon_{n,m}$. The decreasing nature of the sequence $(\epsilon_{n,m})_m$ implies we may assume (with no loss of generality) that
$$E_{n,m}\subset E_{n,m+1}$$
for all $n$ and $m$.
Moreover, since $\mu$ is non-atomic, we may assume $\mu(E_{n,m})=1-1/(2m^n)$ by choosing $E_{n,m}$ as a subset of $E_{\epsilon_{n,m}}$ if necessary. Thus, for every $n$ and each $m\ge 1$,
$$\mu(E_{n,m+1}\setminus E_{n,m})=\tf{1}{2m^n}-\tf{1}{2(m+1)^n}=\tf{(m+1)^n-m^n}{2m^n(m+1)^n}\le \tf{n}{2m^{n+1}}.$$
Let $D_{n,1}=E_{n,1}$ and $D_{n,m}=E_{n,m+1}\setminus E_{n,m}$ for all $m>1$. 

Suppose the above construction has been made for all $n$. 
For each $n$, define a measurable function $g_n$ by
$$g_n=\sum_m1_{D_{n,m}}m^{(n+k)/r}$$
where $k$ is chosen so that $pk/r>1$. Clearly, $g_n\ge 1$ on $\Omega$ for all $n$. 

Let $r_n=\tf{rn}{n+k}$. Then $(r_n)$ increases to $r$. Furthermore, it is easy to see that there is a constant $C<\infty$ which is independent of $n$ so that 
$$\mu(\{g_n>t\})= \sum_{m>t^{r/(n+k)}}\mu(D_{n,m})\le \sum_{m>t^{r/(n+k)}}\tf{n}{2m^{n+1}}\le C t^{-rn/(n+k)}.$$
Thus, $\sup_n\|g_n\|_{r_n,\infty}<\infty$. Moreover, if $f\in \mathcal{F}$ then
$$\mu(\{|f|/g_n>t\})= \sum_m \mu(D_{n,m}\cap\{|f|>tm^{(n+k)/r}\}).$$
Since $D_{n,m}\subset  E_{\epsilon_{n,m}}$, $N3$ and our choice of $k$ imply that there is a constant $C=C(p,q,r)<\infty$ which is independent of $n$ so that
$$\mu(\{|f|/g_n>t\})\le t^{-p}\sum_m 2^{p/r}m^{np/r}m^{-p(n+k)/r}=t^{-p}2^{p/r}\sum_mm^{-pk/r}=Ct^{-p}.$$

Since $r>2$ and $(r_n)$ is increasing to $r$, we may restrict our attention to all $n$ such that $2<r_n<r$. Without loss of generality, we may assume that $r_1>2$. Since $\mu$ is a probability measure 
$$\sup_n\|g_n\|_2\le C(r_1)\sup_n\|g_n\|_{r_n,\infty}<\infty.$$
By the Banach-Saks Theorem for $L_2$, there exists a subsequence $(g_{n_k})_k$ and $g\in L_2$ so that 
$$g=\lim_N G_N\q (in\ L_2)$$
where
$$G_N=N^{-1}\sum_{j=1}^Ng_{n_j}.$$
For any $m\ge 1$ define
$$G_{N,m}=N^{-1}\sum_{j=m}^{N+m-1}g_{n_j}.$$
By a simple application of the triangle inequality, it is easy to see that $g$ is the $L_2$ limit of $(G_{N,m})_N$ for any $m$. By Fatou's lemma,
$$\mu(\{g>t\})\le \liminf_N\mu(\{G_{N,m}>t\})$$
for any $m$.
However, since $(r_n)$ is increasing and $\mu$ is a probability measure, notice that $G_{N,m}$ is an average of functions whose $L_{r_m,\infty}$ norms are uniformly bounded in $m$ by a constant independent of $N$.  So, there exists a constant $C(p,q)<\infty$ so that
$$\sup_m\|g\|_{r_m,\infty}<C(p,q)$$
Since $\mu$ is a probability measure and $(r_m)$ increases to $r$, it follows that $g\in L_{r,\infty}$ with $\|g\|_{r,\infty}\le C(p,q)$.

Fix $f\in \mathcal{F}$. Again, by applying Fatou's lemma,
$$\mu(\{|f|/g>t\})\le \liminf_N\mu(\{|f|/G_N>t\}).$$
By the arithmetic-geometric mean inequality
$$G_N\ge \Big(\prod_{j=1}^{N}g_{n_j}\Big)^{1/N}.$$
Therefore, by H\"older's inequality for weak-type spaces,
$$\mu(\{|f|/G_N>t\})\le t^{-p}\||f|/G_N\|_{p,\infty}^p\le t^{-p}\prod_{i=1}^N\||f|^{1/N}/g_{n_i}^{1/N}\|_{Np,\infty}^p=t^{-p}\prod_{i=1}^N\|f/g_{n_i}\|_{p,\infty}^{p/N}.$$
By normalizing $g$ in $L_{r,\infty}$, the condition $\sup_n\sup_{f\in \mathcal{F}}\|f/g_n\|_{p,\infty}<\infty$ implies there is a constant $C<\infty$ such that
$$\sup_{f\in\mathcal{F}}\|f/g\|_{p,\infty}\le C$$
and $N1$ follows.

\end{proof}



\begin{thebibliography}{10}

\bibitem{AlbiacKalton} F. Albiac and N. Kalton, \emph{Topics in Banach spaces},Springer, New York (2006).

\bibitem{Diestel} G. Diestel, \emph{Factoring multisublinear maps}, J. Functional Analysis {\bf 266} issue 4 (2014), 1928--1947.

\bibitem{Grafakosbook1} L. Grafakos,  \emph{Classical Fourier Analysis,}Second Edition, Graduate Texts in Math., no 249, Springer, New York, 2008.

\bibitem{Kalton2} N.J. Kalton, \emph{Rademacher series and decoupling},New York J. of Math. {\bf 11}(2005), 563--595.

\bibitem{Maurey} B. Maurey, \emph{Th\'eor\`ems de factorisation pour les op\'erateurs lin\'eaires \`a valeurs dans les epaces $L_p$},Soci\'et\'e Mathe\'matique de France, Paris, 1974, With an English summary; Ast\'erisque, No. 11 .(French)

\bibitem{Nikishin} E.M. Nikishin,
\emph{Resonance theorems and superlinear operators}, Uspekhi Mat. Nauk 25 No. 6 (1970), 129--191.(Russian)

\bibitem{Pisier} G. Pisier,
\emph{Factorization of Operators Through $L_{p,\infty}$ or $L_{p,1}$ and Non-Commutative Generalizations,}Math. Ann. {\bf 276} (1986), 105--136.


\bibitem{Wojtaszczyk}P. Wojtaszczyk, \emph{Banach spaces for analysts.}Cambridge Studies in Advanced Mathematics, vol. 25, Cambridge University Press, Cambridge, 1991.

\bibitem{handbook1} \emph{Handbook of the Geometry of Banach Spaces, Vol. 1}. Edited by W.B. Johnson and J. Lindenstrauss. North-Holland Publishing Co., Amsterdam, 2001.

\bibitem{handbook2} \emph{Handbook of the Geometry of Banach Spaces, Vol. 2}. Edited by W.B. Johnson and J. Lindenstrauss. North-Holland Publishing Co., Amsterdam, 2001.
\end{thebibliography}
\end{document}